\documentclass[12pt]{article}
\usepackage{graphicx}
\usepackage{amsmath,amsthm,amssymb,enumerate}%, esint}
\usepackage{euscript,mathrsfs}
\usepackage{color}
\usepackage{dsfont}
\usepackage[left=2cm,right=2cm,top=3.5cm,bottom=3.5cm]{geometry}
\usepackage{color}
\usepackage[framemethod=tikz]{mdframed}
\usepackage{bm}
\allowdisplaybreaks

\usepackage{soul}

\catcode`\@=11 \@addtoreset{equation}{section}

\catcode`\@=12

\newtheorem{Theorem}{Theorem}[section]
\newtheorem{Proposition}[Theorem]{Proposition}
\newtheorem{Lemma}[Theorem]{Lemma}
\newtheorem{Corollary}[Theorem]{Corollary}

\theoremstyle{definition}
\newtheorem{Definition}[Theorem]{Definition}

\newtheorem{Remark}[Theorem]{Remark}

\newcommand{\bTheorem}[1]{
\begin{Theorem} \label{T#1} }
\newcommand{\eT}{\end{Theorem}}

\newcommand{\bProposition}[1]{
\begin{Proposition} \label{P#1}}
\newcommand{\eP}{\end{Proposition}}

\newcommand{\bLemma}[1]{
\begin{Lemma} \label{L#1} }
\newcommand{\eL}{\end{Lemma}}

\newcommand{\bCorollary}[1]{
\begin{Corollary} \label{C#1} }
\newcommand{\eC}{\end{Corollary}}
\newcommand{\bRemark}[1]{
\begin{Remark} \label{R#1} }
\newcommand{\eR}{\end{Remark}}

\newcommand{\bDefinition}[1]{
\begin{Definition} \label{D#1} }
\newcommand{\eD}{\end{Definition}}
\newcommand{\eq}[1]{\begin{equation}
\begin{split}
#1
\end{split}
\end{equation}}

\newcommand{\lr}[1]{\left( #1 \right)}
\newcommand{\Del}{\Delta_x}

\newcommand{\bfphi}{\boldsymbol{\varphi}}

\newcommand{\Td}{\mathbb{T}^d}
\newcommand{\bFormula}[1]{
\begin{equation} \label{#1}}
\newcommand{\eF}{\end{equation}}

\newcommand{\Ov}[1]{\overline{#1}}

\newcommand{\vr}{\varrho}

\newcommand{\vu}{\vc{u}}

\newcommand{\vc}[1]{{\bm #1}}

\newcommand{\Div}{{\rm div}_x}
\newcommand{\Grad}{\nabla_x}

\newcommand{\dx}{\,{\rm d} {x}}

\newcommand{\dt}{\,{\rm d} t }

\newcommand{\intO}[1]{\int_{\Omega} #1 \ \dx}

\newcommand{\intTd}[1]{\int_{\mathbb{T}^d} #1 \ \dx}

\newcommand{\D}{{\rm d}}

\newcommand{\br}{ \nonumber \\ }

\newcommand*\samethanks[1][\value{footnote}]{\footnotemark[#1]}
\def\softd{{\leavevmode\setbox1=\hbox{d}%
          \hbox to 1.05\wd1{d\kern-0.4ex{\char039}\hss}}}
\definecolor{Cgrey}{rgb}{0.85,0.85,0.85}
\definecolor{Cblue}{rgb}{0.50,0.85,0.85}
\definecolor{Cred}{rgb}{1,0,0}
\definecolor{fancy}{rgb}{0.10,0.85,0.10}

\newcommand\Cbox[2]{%
    \newbox\contentbox%
    \newbox\bkgdbox%
    \setbox\contentbox\hbox to \hsize{%
        \vtop{
            \kern\columnsep
            \hbox to \hsize{%
                \kern\columnsep%
                \advance\hsize by -2\columnsep%
                \setlength{\textwidth}{\hsize}%
                \vbox{
                    \parskip=\baselineskip
                    \parindent=0bp
                    #2
                }%
                \kern\columnsep%
            }%
            \kern\columnsep%
        }%
    }%
    \setbox\bkgdbox\vbox{
        \color{#1}
        \hrule width  \wd\contentbox %
               height \ht\contentbox %
               depth  \dp\contentbox
        \color{black}
    }%
    \wd\bkgdbox=0bp%
    \vbox{\hbox to \hsize{\box\bkgdbox\box\contentbox}}%
    \vskip\baselineskip%
}

\mdfdefinestyle{MyFrame}{%
	linecolor=black,
	outerlinewidth=1pt,
	roundcorner=5pt,
	innertopmargin=\baselineskip,
	innerbottommargin=\baselineskip,
	innerrightmargin=10pt,
	innerleftmargin=10pt,
	backgroundcolor=white!20!white}

\makeindex
\begin{document}

%%%%%%%%%%%%%%%%%%%%%%%%%%%%%%%%

\title{Nonuniqueness of weak solutions to the dissipative Aw-Rascle model}

\author{Nilasis Chaudhuri\thanks{The work of N.C. and E.Z. was  supported by the EPSRC Early Career Fellowship no. EP/V000586/1.}
\and Eduard Feireisl
	\thanks{The work of E.F. was partially supported by the
		Czech Sciences Foundation (GA\v CR), Grant Agreement
		21--02411S. The Institute of Mathematics of the Academy of Sciences of
		the Czech Republic is supported by RVO:67985840. }
\and  Ewelina Zatorska\samethanks[1]
}

\date{\today}

\maketitle

	\centerline{$^* $ Department of Mathematics, Imperial College London,   }
	\centerline{ 180 Queen's Gate, SW7 2AZ, London, United Kingdom}
	\vspace{5mm}
\centerline{$^\dagger$  Institute of Mathematics of the Academy of Sciences of the Czech Republic}

\centerline{\v Zitn\' a 25, CZ-115 67 Praha 1, Czech Republic}

\begin{abstract}
	
We prove nonuniqueness of weak solutions to multi-dimensional generalisation of the Aw-Rascle model of vehicular traffic. Our generalisation includes the velocity offset in a form of gradient of density function, which results in a dissipation effect, similar to viscous dissipation in the compressible viscous fluid models. We show that despite this dissipation, the extension of the method of convex integration can be applied to generate infinitely many weak solutions connecting arbitrary initial and final states. We also show that for certain choice of data, ill posedness holds in the class of admissible weak solutions.

\end{abstract}

\bigskip

%\tableofcontents

\section{Motivation} 
\label{i}

The Aw--Rascle model of vehicular traffic (AR) is a second order macroscopic model of traffic developed originally in one-dimensional framework \cite{AR}. It is a system of two conservation laws describing the conservation of mass and conservation of linear momentum. However, unlike in the fluid models, the second equation is not associated with the actual velocity of motion $u$, but with the preferred velocity $w$:
\eq{\label{AR_1D}
		\partial_t \vr + \partial_x(\vr u) &= 0,  \\
		\partial_t (\vr w) + \partial_x (\vr w u ) &= 0, \\ 
		{w} &= {u} + P(\vr) . 
}	
The two velocities differ by the velocity offset denoted by $P(\vr)>0$, which, in this particular case, depends only on the density. The relation $u=w-P(\vr)$ means that the actual velocity of motion is always smaller 	than the preferred velocity depending on the congestion of the cars ahead.
The one-dimensional AR system has been derived in \cite{AwKlar} from the particle model called Follow-the-Leader model with particular form of the offset function $P(\vr)=\vr^\gamma$.
% (FtL) particle model
% \eq{\label{FtL}
% 		\dot x_i &= u_i,\\ 
% 				\dot u_i &= C\frac{u_{i+1}-u_i}{(x_{i+1}-x_i)^{\gamma+1}}, 
% 		}
% describing the evolution of positions $x_i$ and velocities $u_i$ of the cars on the single lane road for an ensemble of cars $i=1,\ldots,N$. The exponent in the denominator of the second equation corresponds to the form of the offset function $P(\vr)$ in the macroscopic system, and in this particular case it is equal to $\vr^\gamma$.
This form of the offset function has several drawbacks. First of all, the maximal velocity and maximal density constraint are not preserved by the system \eqref{AR_1D}. Secondly, it is not clear what this offset should be for multi-dimensional model of traffic, where the velocity is a vector and the offset is a scalar function.

The remedy to the first problem has been proposed \cite{BDDR} where the authors considered a cost function with maximal density constraint $\bar\vr>0$
\eq{\label{sing_cost}
P(\vr)=\lr{\frac{1}{\vr}-\frac{1}{\bar\vr}}^{-\gamma}.}
In this way, the density $\vr$ stays always below its critical value $\bar\vr$, provided it was so initially.
% In \cite{BDDR} it was also pointed out, that the velocity offset $P(\vr)$ should be very small unless the density $\vr$ is very close to the maximal value $\bar\vr$. It was advocated by the fact that the drivers do not reduce their speed significantly if the traffic is not congested enough. To investigate this phenomenon they considered the singular limit of the model \eqref{AR_1D} with $P(\vr)$ replaced by $\ep P(\vr)$ and $\ep\to 0$.\\
%, which (formally) gives rise to the constrained pressureless gas dynamics model.\\
To talk about the remedy to the second problem -- the dimension discrepancy -- we first point out that avoidance of collisions can also be modelled by introduction of the force that becomes singular at the contact points. This is the main idea behind the macroscopic models of lubrication, considered for example in \cite{LM}. In particular, the one-dimensional macroscopic lubrication model for interacting rigid spheres integrating the inertial effects reads:
\eq{\label{NS_1D}
&\partial_t \vr + \partial_x(\vr u) = 0,\\
&\partial_t(\vr u)+\partial_x(\vr u^2)-\partial_x\lr{\mu(\vr)\partial_x u}=0.}
Here $\mu(\vr)>0$  is a singular function of the density such that  $\mu(\vr)\to \infty$ when density tends to some maximal constraint, say $\bar\vr$, which is dictated by the physical dimensions of the interacting balls. Note that \eqref{NS_1D} is in fact the compressible, pressureless Navier-Stokes system, studied for example in \cite{HZ}. Extension of this system to multi-dimensional case, in particular, the form of the stress and maximal constraint of the density are not known. Note, however, that 
% taking in this system
% \eq{
% p'(\vr)=\frac{\mu(\vr)}{\vr^2}
% }
a simple formal calculation allows to convert the system \eqref{NS_1D} into a version AR system \eqref{AR_1D} with 
\eq{P(\vr)=\frac{\mu(\vr)}{\vr^2}\partial_x \vr.}
Taking this form of the offset function has the advantage of dimension compatibility with the velocity vector field in higher dimensions. The existence of measure-valued solutions to such generalisation in multi-dimensional setting, and their weak-strong uniqueness was recently proved in \cite{CGZ}.  

Replacement of a scalar function by its gradient in the form of offset $P(\vr)$ accounts for including certain non-local effects in the interactions between the drivers at the microscopic level in the Follow-the-Leader model.
% . Indeed, to preserve the constraint $\dot{w}_i=0$ in the FtL particle model \eqref{FtL}, one would have to include some effect on the acceleration of $i$-th vehicle by the vehicles $i+1$ and $i+2$. 
Another generalisation of the AR model including the non-local interactions was recently studied in \cite{CFGG}.

On the other hand, in line with the derivation from \cite{HMV}, the velocity offset in the two-dimensional AR system is a vector of functions i.e. $P(\vr)=\vc{h}(\vr)=[h_1(\vr),h_2(\vr)]$.

In the current paper we consider a combination of these ansatz and we study a $d$-dimensional AR model of the form:
\begin{mdframed}[style=MyFrame]
	\begin{align}
		\partial_t \vr + \Div (\vr \vu) &= 0, \label{i1} \\
		\partial_t (\vr \vc{w}) + \Div (\vr \vc{w} \otimes \vu ) &= 0, \label{i2} \\ 
		\vc{w} &= \vc{u} + \vc{h}(\vr) + \Grad p(\vr) . \label{i3}
	\end{align}	
		\end{mdframed}

For simplicity, we consider the periodic boundary conditions -- the physical domain is identified with 
the $d-$dimensional flat torus 
\begin{equation} \label{i4}
	\Td = \left( [-1,1]|_{\{ -1; 1 \} } \right)^d,\ d=2,3.
	\end{equation}

	The goal of the paper is to show that, similarly to the compressible Euler system, the Aw--Rascle system is basically ill--posed in the class of weak (distributional) solutions.
To this end, we adapt the general approach developed in \cite{Fei2016} based on the method of {\emph{convex integration}}. 
This method was introduced by DeLellis and Sz\'ekelyhidi \cite{DeSz}, primarily to prove the existence of infinitely many {\emph{wild}} solutions to the incompressible Euler system. Subsequently, it was extended by Chiodaroli \cite{Ch} to the compressible Euler system, and more recently by  Buckmaster and Vicol for the incompressible Navier-Stokes equations \cite{BuVi}. It is not yet known if convex integration technique could be further extended to the compressible Navier-Stokes equations. Note, however, that weak inviscid limit of  compressible Navier--Stokes system with degenerate viscosities has been recently used in \cite{CVY} to generate infinitely many global-in-time admissible weak solutions to the isentropic Euler system.
The fact that the convex integration technique works for system (\ref{i1})--(\ref{i3}), which in one-dimensional setting coincides with the compressible Navier-Stokes system \eqref{NS_1D}, is therefore a very interesting observation. 
In multi-dimensional setting system (\ref{i1})--(\ref{i3}) is equivalent to a dissipate pressureless compressible system with degenerate, density-dependent shear viscosity and a lower order drift term. Indeed, taking $ \vc{h}=0 $  in (\ref{i3}), and substituting for $\vc{w}$ in the equation \eqref{i2}, we formally obtain 
\begin{align*}
    % \partial_t \vr + \Div (\vr \vu) &= 0,  \\
		\partial_t (\vr \vc{u}) + \Div (\vr \vc{u} \otimes \vu ) = \nabla_x (\vr Q^\prime (\vr)\Div \vu) + \mathcal{L}[\nabla_x Q(\vr), \nabla_x \vu] ,
\end{align*}
where $ Q^\prime(\vr) = \vr p^\prime (\vr)$ and $$\mathcal{L}[\nabla_x Q(\vr), \nabla_x \vu] = \nabla_x(\nabla_x Q(\vr) \cdot \vu )- \Div(\nabla_x Q(\vr) \otimes \vu ),$$ 
which is a lower order term, a simple calculation yields
\begin{align*}
     \left(\mathcal{L}[\nabla_x Q(\vr), \nabla_x \vu]\right)_j= \sum_{i=1}^{3} \left( \partial_{x_i}Q(\vr) \partial_{x_j}u_{i} - \partial_{x_j} Q(\vr) \partial_{x_i} u_i   \right)  \text{ for } j=1,2,3.
\end{align*}
{We can thus say that our idea  works for certain viscous compressible models with degenerate viscosity coefficients possessing the ``two-velocity" structure''. Similar structure has been used in the past to prove the existence of solutions to compressible Navier-Stokes equations with density-dependent viscosity \cite{BVY}, and in \cite{BFHZ} to consider stochastically perturbed transport terms in the compressible Navier-Stokes system with constant viscosity coefficients.}

	The paper is organised as follows. In Section \ref{Sec:2} we state our first main result, Theorem \ref{ccT1}, about ill posedness of the Aw--Rascle system \eqref{i1} -- \eqref{i3} with respect to the initial-final data. The solutions obtained in this section connect \emph{arbitrary} initial and terminal states, however, 
they may violate the energy inequality. The ill posedness in the class of {\emph{admissible}} weak solutions satisfying this inequality is shown in Section \ref{Sec:3}, the final result is stated in  Theorem \ref{ccT2}. The paper is concluded with a discussion of other boundary conditions.

% \section{Generalized Aw-Rascle system, convex integration}
% \label{c}

% The goal of the paper is to show that, similarly to the compressible Euler system, the Aw--Rascle system is basically ill--posed in the class of weak (distributional) solutions.
% To this end, we adapt the general approach developed in \cite{Fei2016} based on the method of convex integration. 

\section{Ill posedness with respect to the initial--final data}\label{Sec:2}
In this section we formulate and prove our  first main result: that any initial density--velocity data $(\vr_0, \vu_0) = (\vr(0, \cdot), \vu(0, \cdot))$ 
can connect to arbitrary terminal state $(\vr_T, \vu_T) = ((\vr(T, \cdot), \vu(T, \cdot))$ via a weak solution 
to problem \eqref{i1}--\eqref{i4}. More specifically, we consider 
\begin{equation} \label{cc1}
	\vr_0, \vr_T \in C^2(\Td),\ \inf_{\Td} \vr_0 > 0,\ \inf_{\Td} \vr_T > 0, \ \intTd{ \vr_0 } = \intTd{ \vr_T }
	\end{equation}
together with 
\begin{equation} \label{cc2}
	\vu_0, \vu_T \in C^3(\Td; R^d),\ \intTd{ \vr_T \vu_T } - \intTd{\vr_0 \vu_0 } = \intTd{ \vr_0 \vc{h} (\vr_0) } - 
	\intTd{ \vr_T \vc{h}(\vr_T) }.
	\end{equation}
Note that the integral equalities in \eqref{cc1}, \eqref{cc2} represent necessary compatibility conditions as the quantities 
\[
\intTd{ \vr(t, \cdot) },\ \intTd{ \vr \vc{w} (t, \cdot)}
\]
are conserved even in the class of weak solutions. 

We claim the following result.

	\begin{Theorem}[{\bf Ill posedness with respect to the data}] \label{ccT1}
		Let $d=2,3$.
		Suppose that 
		\begin{equation} \label{hyp1}
		\vc{h} \in C^2 ((0, \infty) ; R^d),\ p \in C^2((0,\infty)).
		\end{equation}
		Let $(\vr_0, \vu_0)$, $(\vr_T, \vu_T)$ satisfy \eqref{cc1}, \eqref{cc2}. 
		
		Then the system \eqref{i1} -- \eqref{i3}, endowed with the periodic boundary  conditions 
		\eqref{i4} admits infinitely many weak solutions in the class 
		\[
		\vr \in C^2([0,T] \times \Td), \vu \in L^\infty ((0,T) \times \Td; R^d)
		\]
		such that
		\begin{equation} \label{cc3}
			\vr(0, \cdot) = \vr_0,\ \vr(T, \cdot) = \vr_T,\ 
			(\vr \vu) (0, \cdot) = \vr_0 \vu_0,\ 
			(\vr \vu) (T, \cdot) = \vr_T \vu_T.
		\end{equation}

		\end{Theorem}

\begin{Remark} \label{rema1}
	
	It will become clear in the course of the proof 
		(see formula \eqref{conf} below) that hypothesis \eqref{hyp1} can be relaxed to
	\[
	\vc{h} \in C^2(I; R^d),\ p \in C^2(I; R^d),
	\]  
where $I \subset (0, \infty)$ is an open interval containing the convex closure of the range 
		of $\vr_0$, $\vr_T$. In particular, we can choose  $\vc{h}(\vr),p(\vr)$ to be singular as in the form \eqref{sing_cost} proposed by the authors of \cite{BDDR}.
	
	\end{Remark}

Here and hereafter, we adopt the standard definition of \emph{weak solution} via the integral identities: 
\begin{align} 
	\int_0^T \intTd{ \left[ \vr \partial_t \varphi + \vr \vu \cdot \Grad \varphi \right] } \dt &= 0 \br 
\mbox{for any} \ \varphi &\in C^1_c((0,T) \times \Td); \br
\int_0^T \intTd{\left[  \vr \vc{w} \cdot \partial_t \bfphi + \vr \vc{w} \otimes \vu : \Grad \bfphi \right] } \dt& = 0 \br 
\mbox{for any}\ \bfphi &\in C^1_c((0,T) \times \Td; R^d).
\label{cc4}
\end{align}
In particular, 
\[
\vr \vc{w} \in C_{\rm weak}([0,T]; L^q(\Td; R^d)) \ \mbox{for any}\ 1 \leq q < \infty; 
\]
hence 
\[
\vr \vc{u} \in C_{\rm weak}([0,T]; L^q(\Td; R^d)) \ \mbox{for any}\ 1 \leq q < \infty, 
\]
and \eqref{cc3} makes sense. 

In the remaining part of the paper, we develop an abstract framework that enables to prove Theorem \ref{ccT1} along with other results stated below.

\subsection{Momentum decomposition}

Write 
\begin{equation} \label{c1}
	\vr \vu = \vc{v} + \vc{V} + \Grad \Phi , 
	\end{equation}
where 
\begin{equation} \label{c2}
	\Div \vc{v} = 0,\ \intTd{ \vc{v} } = 0,\ 
	\vc{V} = \vc{V}(t) \in R^d.
\end{equation}
Accordingly, the equation of continuity \eqref{i1} reads 
\begin{equation} \label{c3}
	\partial_t \vr + \Del \Phi = 0.
\end{equation}

\subsubsection{Density profile}
\label{rem1}

The next step is adjusting a suitable density profile, 
\begin{equation} \label{c4}
	\vr \in C^2([0,T] \times \Td),\ \vr(0, \cdot) = \vr_0, \vr > 0, \vr(T, \cdot) = \vr_T
\end{equation}
where $\vr_0$, $\vr_T$ are the desired initial and terminal states.
In accordance with \eqref{c3}, this should be done in such a way that
\begin{equation} \label{c5}
	\partial_t \vr (0, \cdot) + \Del \Phi_{0} = 0,\ 	\partial_t \vr (T, \cdot) + \Del \Phi_{T} = 0,
\end{equation}
where $\Phi_0$, $\Phi_T$ are the values of the acoustic potential determined by the Helmholtz decomposition of the initial data, and terminal data
\begin{equation} \label{c6}
	\vr_0 \vu_0 = \vc{v}_0 + \vc{V}_0 + \Grad \Phi_0, \ \vr_T \vu_T = \vc{v}_T + \vc{V}_T + \Grad \Phi_T,
\end{equation}
respectively.

Consider the functions 
\begin{align}
H &\in C^\infty[0,T],\ 0 \leq H \leq 1,  H (0) = 1, H (T) = 0, H' (0) = H' (T) = 0, \br
Z_0^\delta &\in C^\infty_c [0,T),\ Z_0^\delta (0) = 0,\ (Z_0^\delta)'(0) = -1, \br   
Z_T^\delta &\in C^\infty_c (0,T],\ Z_T^\delta (T) = 0,\ (Z_T^\delta)'(T) = -1, \br
|Z_0^\delta|,\ |Z_T^\delta| &< \delta, \ \delta > 0.    
\nonumber
\end{align}

The desired density profile can be taken as 
\begin{equation} \label{conf}
\vr(t,x) = H(t) \vr_0 (x) + \vr_T (x) (1 - H(t) ) + Z_0^\delta(t) \Del \Phi_0(x) + Z_T^\delta (t) \Del \Phi_T(x).
\end{equation}
Indeed it is easy to check that \eqref{c4}, \eqref{c5} hold while the acoustic potential $\Phi$ is uniquely determined by \eqref{c3}.
Moreover, if $\delta > 0$ is chosen small enough, we get 
\begin{equation} \label{c6b}
	\inf_{(0,T) \times \Td} \vr > 0. 
	\end{equation}

\subsection{Transformed problem I}
	
With $\vr$, $\Phi$ fixed in the preceding part, the problem \eqref{i1}--\eqref{i3} reduces to 

% \begin{mdframed}[style=MyFrame]
	
\begin{align}
\partial_t (\vc{v} + \vc{V}) &+ \Div \left( \frac{ \left( \vc{v} + \vc{V} + \Grad \Phi  
	 \right) \otimes  \left( \vc{v} + \vc{V} + \Grad \Phi \right) }{\vr} + \partial_t \left( \Phi + 
 P(\vr) \right) \mathbb{I} \right) \br 	&= - \partial_t (\vr \vc{h}(\vr))  - \Div \Big( \left(\vc{h}(\vr) + \Grad p(\vr) \right) \otimes \left( \vc{V} + \Grad \Phi \right) \Big) \br 
 &\quad - \Grad (\vc{h}(\vr) + \Grad p(\vr) ) \cdot \vc{v} 
\label{c8} \\ 
\Div \vc{v} &= 0,
\label{c9}
	\end{align}

% \end{mdframed}

\noindent where $\Grad P(\vr) = \vr \Grad p(\vr)$.
System \eqref{c8}, \eqref{c9} still contains two unknowns -- $\vc{v}$ and $\vc{V}$ that should satisfy the associated initial and terminal conditions 
\[
\vc{v}(0, \cdot) = \vc{v}_0,\ \vc{V}(0, \cdot) = \vc{V}_0,\ 
\vc{v}(T, \cdot) = \vc{v}_T,\ \vc{V}(T, \cdot) = \vc{V}_T.
\]

\subsection{Fixing $\vc{V}$}

In addition to \eqref{c5}, the density profile should give rise to the desired momentum average $\vc{V}$. In accordance with the momentum equation 
\eqref{i2}, we get 
\begin{equation} \label{c10}
	\vc{V}(t) = \vc{V}_0 -  \frac{1}{|\Td|} \int_0^t \intTd{ \partial_t (\vr \vc{h}(\vr)) (s, \cdot) } \D s
\end{equation}
so that 
\[
\partial_t \vc{V} = - \frac{1}{|\Td|} \intTd{ \partial_t (\vr \vc{h}(\vr)) (t, \cdot) }.
\]
Moreover, in accordance with \eqref{cc2}, \eqref{cc3}, we have 
\begin{align}
	\intTd{ \vr_T \vu_T }  &=|\Td |	\vc{V}_T = |\Td | \vc{V}(T) = |\Td | \vc{V}_0 -  \int_0^T \intTd{ \partial_t (\vr \vc{h}(\vr)) (s, \cdot) } \D s \br
	&= \intTd{ \vr_0 \vu_0 } - \intTd{ \vr_T \vc{h} (\vr_T) } + \intTd{ \vr_0 \vc{h}(\vr_0)}.
	\label{cc5}
\end{align}	

Consequently, 
equation \eqref{c8} reduces to
\begin{align}
	\partial_t \vc{v} &+ \Div \left( \frac{ \left( \vc{v} + \vc{V} + \Grad \Phi  
		\right) \otimes  \left( \vc{v} + \vc{V} + \Grad \Phi \right) }{\vr} + \partial_t \left( \Phi + 
	P(\vr) \right) \mathbb{I} \right) \br 	&= \left( \frac{1}{|\Td|} \intTd{ \partial_t (\vr \vc{h}(\vr)) }  - \partial_t (\vr \vc{h}(\vr))  \right) - \Div \Big( \left(\vc{h}(\vr) + \Grad p(\vr) \right) \otimes \left( \vc{V} + \Grad \Phi \right) \Big) \br 
	&\quad - \Grad (\vc{h}(\vr) + \Grad p(\vr) ) \cdot \vc{v} 
	\label{c11} 
	\end{align}
Similarly to the acoustic potential $\Phi$, the function $\vc{V}$ is now determined through the given density profile $\vr$ via \eqref{c10}.

\subsection{Elliptic problem I}

To rewrite \eqref{c11} in the form considered in \cite{Fei2016}, we consider
a symmetric traceless tensor 
\begin{equation} \label{c12}
\mathbb{M} = \Grad \vc{U} + \Grad \vc{U}^t - \frac{2}{d} \Div \vc{U} \mathbb{I}, 
\end{equation}
where $\vc{U}$ is the unique zero--mean solution of the elliptic problem 
\begin{align}
\Div \left( \Grad \vc{U} + \Grad \vc{U}^t - \frac{2}{d} \Div \vc{U} \mathbb{I}	\right) &= 
\Div \Big( \left(\vc{h}(\vr) + \Grad p(\vr) \right) \otimes \left( \vc{V} + \Grad \Phi \right) \Big) \br 
&-\left( \frac{1}{|\Td|} \intTd{ \partial_t (\vr \vc{h}(\vr)) }  - \partial_t (\vr \vc{h}(\vr))  \right).
\label{c13}
\end{align}
Consequently, problem \eqref{c8}, \eqref{c9} can be rewritten in the form

% \begin{mdframed}[style=MyFrame]
	
	\begin{align}
		\partial_t \vc{v} &+ \Div \left( \frac{ \left( \vc{v} + \vc{V} + \Grad \Phi  
			\right) \otimes  \left( \vc{v} + \vc{V} + \Grad \Phi \right) }{\vr} + \partial_t \left( \Phi + 
		P(\vr) \right) \mathbb{I} + \mathbb{M} \right) \br  
		&= - \Grad (\vc{h}(\vr) + \Grad p(\vr) ) \cdot \vc{v}, 
		\label{c14} \\ 
		\Div \vc{v} &= 0, \label{c15a}
	\end{align}
	with $\vc{v}(0, \cdot) = \vc{v}_0$, $\vc{v}(T, \cdot) = \vc{v}_T$.

% \end{mdframed}

\subsection{Elliptic problem II}

Similarly to the preceding step, we set 
\begin{equation} \label{c16}
	\mathbb{N} = \Grad \vc{R} + \Grad \vc{R}^t - \frac{2}{d} \Div \vc{R} \mathbb{I}, 
\end{equation}
with $\vc{R}$ solving 
\begin{equation} \label{c17}
	\Div \left( \Grad \vc{R} + \Grad \vc{R}^t - \frac{2}{d} \Div \vc{R} \mathbb{I} \right) = 
\Grad (\vc{h}(\vr) + \Grad p(\vr) ) \cdot \vc{v}.
\end{equation}	
Note carefully that $\mathbb{N} = \mathbb{N}[\vc{v}]$ depends on the unknown $\vc{v}$.

Consequently, we may rewrite \eqref{c14}, \eqref{c15a} in the form
	
	\begin{align}
		\partial_t \vc{v} + \Div \left( \frac{ \left( \vc{v} + \vc{V} + \Grad \Phi  
			\right) \otimes  \left( \vc{v} + \vc{V} + \Grad \Phi \right) }{\vr} + \partial_t \left( \Phi + 
		P(\vr) \right) \mathbb{I} + \mathbb{M} + \mathbb{N}[\vc{v}] \right)  &= 0, 		
		\label{c18} \\ 
		\Div \vc{v} &= 0.
		\label{c19}
	\end{align}

\subsection{Adjusting the energy}

Finally, let us consider the energy associated to the system, 
\begin{equation} \label{c20}
e = \frac{1}{2} \frac{|\vc{v} + \vc{V} + \Grad \Phi |^2 }{\vr}.	
	\end{equation}
Introducing the notation 
\[
\vc{m} \odot \vc{m} = \vc{m} \otimes \vc{m} - \frac{1}{d} |\vc{m}|^2 \mathbb{I},
\]
we may rewrite \eqref{c18}, \eqref{c19} as an abstract ``Euler system'':

% \begin{mdframed}[style=MyFrame]
	\begin{align}
		\partial_t \vc{v} + \Div \left( \frac{ \left( \vc{v} + \vc{V} + \Grad \Phi  
			\right) \odot  \left( \vc{v} + \vc{V} + \Grad \Phi \right) }{\vr} + \mathbb{M} + \mathbb{N}[\vc{v}] \right)  &= 0, 		
		\label{c21} \\ 
		\Div \vc{v} &= 0,
		\label{c22} \\
\frac{1}{2} \frac{|\vc{v} + \vc{V} + \Grad \Phi |^2 }{\vr} = e &= \Lambda - \frac{d}{2} \partial_t 
\Big( \Phi + P(\vr) \Big),		
		\label{c23}\\ 
		\vc{v}(0, \cdot) = \vc{v}_0,\ \vc{v}(T, \cdot) = \vc{v}_T
		\label{cc24}
	\end{align}
	
% \end{mdframed}

\noindent where $\Lambda = \Lambda(t)$ is an arbitrary spatially homogeneous function to be adjusted below. 

\subsection{Convex integration}

Motivated by \cite[Section 13.2.2]{Fei2016}, we introduce the class of \emph{subsolutions} to problem \eqref{c21}--
\eqref{cc24}:

\begin{align}
	X_0 &= \left\{ \vc{v} \in C_{\rm weak}([0,T]; L^2(\Td; R^d) \cap L^\infty((0,T) \times \Td; R^d ) \ \Big|\ 	\vc{v}(0, \cdot) = \vc{v}_0,\ \vc{v}(T, \cdot) = \vc{v}_T, \right. \br
	&\partial_t \vc{v} + \Div \mathbb{F} = 0 \ \mbox{in}\ \mathcal{D}'((0,T) \times \Td; R^d)
	\ \mbox{for some}\ \mathbb{F} \in L^\infty((0,T) \times \Td ; R^{d \times d}_{0,{\rm sym}} ),  \br
	&\vc{v} \in C((0,T) \times \Td; R^d),\ \mathbb{F} \in C((0,T) \times \Td ; R^{d \times d}_{0,{\rm sym}} ), \br
	& \sup_{\tau < t \leq T, \ x \in \Td } \frac{d}{2} \lambda_{\rm max} 
	\left[ \frac{(\vc{v} + \vc{V} + \Grad \Phi) \otimes (\vc{v} + \vc{V} + \Grad \Phi)}{\vr} - 
	\mathbb{F} + \mathbb{M} + \mathbb{N}[\vc{v}] \right] - e < 0 \br  &\mbox{for any} \ 0 < \tau < T \Big\}.
	\label{cc25}	
	\end{align}
In accordance with \eqref{c23}, the energy $e$ in \eqref{cc25} is given as
\[
e = \Lambda - \frac{d}{2} \partial_t \left( \Phi + P(\vr) \right).
\]
The symbol $\lambda_{\rm max}[\mathbb{A}]$ denotes the maximal eigenvalue of a symmetric matrix $\mathbb{A}$. We recall the algebraic inequality
\begin{equation} \label{cc26}
\frac{1}{2} |\vc{w}|^2 \leq d \lambda_{\rm max} [\vc{w} \otimes \vc{w} - \mathbb{B} ],\ 
\mathbb{B} \in R^{d \times d}_{0, {\rm sym}} .
\end{equation}

As proved in \cite[Theorem 13.2.1]{Fei2016}, problem \eqref{c21}--\eqref{cc24} admits infinitely many weak solution
if the following holds:

\begin{itemize}
	
	\item the set $X_0$ of subsolutions is non--empty;
	\item the set $X_0$ is bounded in $L^\infty((0,T) \times \Td; R^d)$; 
	\item the mapping 
	\[
	\vc{v} \mapsto \mathbb{N}[\vc{v}] 
	\]
	enjoys the following weak continuity property: 
	\begin{align}
	\vc{v}_n \to \vc{v} \ \mbox{in}\ C_{\rm weak}([0,\tau]; L^2(\Td; R^d)) \ 
	&\mbox{and weakly-(*) in}\ L^\infty ((0,\tau) \times \Td; R^d)) \br 
	&\Rightarrow \br 
	\mathbb{N}[\vc{v}_n] \to \mathbb{N}[\vc{v}] \ &\mbox{in}\ C([0,\tau] \times \Td ; R^{d \times d})
	\label{cc27}
	\end{align}
	for any $0< \tau \leq T$.
	
	\end{itemize}

To see that $X_0$ is non--empty, it is enough to consider 
\[
\vc{v} = (1 - t/T) \vc{v}_0 + t/T \vc{v}_T 
\]
the obviously satisfies the initial--terminal conditions, $\Div \vc{v} = 0$, and 
\[
\partial_t \vc{v} = \frac{1}{T} (\vc{v}_T - \vc{v}_0 ).
\]
Since 
\[
\intTd{ (\vc{v}_T - \vc{v}_0 ) } = 0, 
\]
it is easy to find (smooth) $\mathbb{F} \in L^\infty((0,T) \times \Td; R^{d \times d}_{0, {\rm sym}})$ such that
\[
\partial_t \vc{v} + \Div \mathbb{F} = 0.
\]

Finally, we fix $\Lambda > 0$ large enough yielding $\vc{v} \in X_0$ -- the set of subsolutions is non--empty. 
Moreover, with $\Lambda$ fixed, we may use inequality \eqref{cc26} to concluded that $X_0$ is bounded 
in $L^\infty((0,T) \times \Td; R^d)$. The continuity property \eqref{cc27} follows easily from \eqref{c16}, \eqref{c17} and the standard elliptic $L^p$--theory. 

We have proved Theorem \ref{ccT1}. $\Box$

\section{Satisfaction of the energy inequality}\label{Sec:3}

The AR system \eqref{i1} -- \eqref{i3} admits a natural energy functional 
\[
E(\vr, \vu) = \frac{1}{2} \vr \left| \vu + \vc{h}(\vr) + \Grad p(\vr) \right|^2. 
\] 
Given the periodic boundary conditions, the total energy of \emph{smooth} solutions is conserved, 
\[
\intTd{ E(\vr, \vu)(t, \cdot) } =  \intTd{ E(\vr_0, \vu_0) } \ \mbox{for any}\ t \in [0,T].
\]
Admissible \emph{weak} solutions should satisfy at least  the energy inequality 
\begin{equation} \label{cc31}
	\frac{\D }{\dt} \intTd{ E(\vr, \vu) } \leq 0,\quad \intTd{ E(\vr, \vu) (t, \cdot) } \leq 
	\intTd{ E(\vr_0, \vu_0 ) }.
\end{equation}	

% \textcolor{red}{Ed: You may prefer the weak formulation}
% \[
% - \int_0^T \intTd{ E(\vr, \vu) \partial_t \psi } \dt \leq \intTd{ E(\vr_0, \vu_0 ) } 
% \]
% \textcolor{red}{for any $\psi \in C^1_c[0,T)$, $\psi \geq 0$, $\psi(0) = 1$.}

\begin{Remark}
The first inequality in \eqref{cc31} is satisfied in the sense of distributions, and the second one guarantees that there is no initial energy jump. Equivalently, we can include both inequalities in a single weak formulation 
\[
- \int_0^T \intTd{ E(\vr, \vu) \partial_t \psi } \dt \leq \intTd{ E(\vr_0, \vu_0 ) } 
\]
satisfied for any $\psi \in C^1_c[0,T)$, $\psi \geq 0$, $\psi(0) = 1$.
\end{Remark}

The solutions obtained in Theorem \ref{ccT1} connect \emph{arbitrary} initial and terminal states, in particular, 
they may violate at least one of the inequalities in \eqref{cc31}.

To obtain the existence of infinitely many admissible solutions, we change slightly the ansatz in Theorem 
\eqref{ccT1} choosing 
\[
\vr_0 = \vr_T,\ \vu_0 = \vu_T = 0. 
\]
Keeping the notation of Section 
\ref{Sec:2} we therefore obtain 
\[
\vr \vu = \vc{v}, \ \Phi = 0,\ \vc{V} = 0, 
\]
while system \eqref{c21}--\eqref{c23} reduces to
% \begin{mdframed}[style=MyFrame]
	
	\begin{align}
		\partial_t \vc{v} + \Div \left( \frac{ \vc{v}  
			 \odot  \vc{v}   }{\vr} + \mathbb{M} + \mathbb{N}[\vc{v}] \right)  &= 0, 		
		\label{c24} \\ 
		\Div \vc{v} &= 0,
		\label{c25} \\
		\frac{1}{2} \frac{|\vc{v} |^2 }{\vr} = e &= \Lambda.		
		\label{c26}
	\end{align}
	
% \end{mdframed}

Seeing that $\vc{v} = \vr \vu$ with $\vr$ independent of time, we have to fix $\Lambda$ in \eqref{c26} so that 
\begin{equation} \label{c27}
	\frac{\D}{\dt} \intTd{ \frac{1}{2} \vr \left|\vu + \vc{h}(\vr) + \Grad p(\vr) \right|^2 } \leq 0.
\end{equation}
We have 
\[
\frac{1}{2} \vr \left|\vu + \vc{h}(\vr) + \Grad p(\vr) \right|^2 = 
\frac{1}{2} \frac {|\vc{v}|^2}{\vr} + \vr \vu \cdot \vc{h}(\vr) + \vr \vu \cdot \Grad \vr + 
\frac{1}{2} \vr |\vc{h}(\vr) + \Grad p(\vr) |^2. 
\]
As $\vr = \vr_0(x)$ is independent of $t$, we easily compute 
\[
\frac{\D }{\dt} \intTd{ E(\vr, \vu) } = \frac{|\Td|}{2} \Lambda'(t) + \frac{\D }{\dt} \intTd{\vr \vu \cdot \vc{h}(\vr)}, 
	\]
where we have used	
\[
\intTd{ \vr \vu \cdot \Grad p(\vr) } = 0.
\]
Finally, using the momentum equation \eqref{i2} we compute 
\[
 \frac{\D }{\dt} \intTd{\vr \vu \cdot \vc{h}(\vr)} = \intTd{ \vr (\vu + \vc{h}(\vr) + \Grad p(\vr) ) \otimes 
 \vu : \Grad \vc{h}(\vr) }.	
\]
Consequently, we may fix $\Lambda = \Lambda (t)$ in such a way that 
\[
\frac{\D }{\dt} \intTd{ E(\vr, \vu) } \leq 0. 
\]

Thus, fixing $\Lambda$ and applying Theorem \ref{ccT1}, we obtain infinitely many solutions of the Aw--Rascle system with a non--increasing total energy profile. This, however, does not exclude the possibility that the 
energy experiences initial jump, specifically, 
\[
\liminf_{t \to 0 +} \intTd{ E(\vr(t), \vu(t)) } > \intTd{ E(\vr_0, \vu_0) }.
\]

To solve the problem of initial energy jump, we use \cite[Theorem 13.6.1]{Fei2016}. Specifically, there exists a sequence of times $\tau_n \searrow 0$ such that the problem \eqref{c24}--\eqref{c26} admits infinitely many weak 
solutions on the interval $[\tau_n, T]$, with the initial data 
\[
(\vr(\tau_n, \cdot), \vu (\tau_n, \cdot)) = (\vr_0, \vu(\tau_n, \cdot))
\]
such that 
\[
\liminf_{t \to \tau_n +} \intTd{ E(\vr(t), \vu(t)) } > \intTd{ E(\vr(\tau_n, \cdot), \vu(\tau_n, \cdot) ) }.
\]

We have shown the following result.

	\begin{Theorem}[{\bf Ill posedness in the class of admissible solutions}] \label{ccT2}
		Let $d=2,3$.
		Suppose that 
		\begin{equation} \label{hyp2}
		\vc{h} \in C^2 ((0, \infty) ; R^d),\ p \in C^2((0,\infty)).
		\end{equation}
		Let $\vr_0 \in C^2(\Td)$, $\inf_{\Td} \vr_0 > 0$ be given. 
		Then there exists an initial velocity $\vu_0 \in L^\infty(\Td; R^d)$  such that
		the system \eqref{i1} -- \eqref{i3}, endowed with the periodic boundary  conditions 
		\eqref{i4} admits infinitely many weak solutions in the class 
		\[
		\vr \in C^2([0,T] \times \Td), \vu \in L^\infty ((0,T) \times \Td; R^d)
		\]
		satisfying
		\begin{equation} \label{cc3b}
			\vr(0, \cdot) = \vr(T, \cdot) = \vr_0,\  
			(\vr \vu) (T, \cdot) = 0,
		\end{equation}
		together with the energy inequality 		
\[		
	\frac{\D }{\dt} \intTd{ E(\vr, \vu) } \leq 0,\ \intTd{ E(\vr, \vu) (t, \cdot) } \leq 
\intTd{ E(\vr_0, \vu_0 ) }.
\]	
% \textcolor{red}{Ed: The same remark as above}	
	\end{Theorem}

\begin{Remark} \label{rema2}
	
Hypothesis \eqref{hyp2} can be relaxed exactly as in Remark \ref{rema1}.

	\end{Remark}

% \section{Other boundary conditions}
	
We conclude the paper with two remarks concerning other choices of boundary data.
\begin{Remark}
The periodic boundary data can be replaced by  more physically realistic boundary conditions, namely 
\begin{equation} \label{I1}
	\vu \cdot \vc{n}|_{\partial \Omega} = 0, 
	\end{equation}
where $\Omega \subset R^d$ is a bounded domain with smooth boundary. Accordingly, the weak formulation of the equation \eqref{i2} reads
\begin{equation} \label{I2}
	\int_0^T \intO{ \Big( \vr \vc{w} \cdot \partial_t \bfphi + 
		\vr \vc{w} \otimes \vu : \Grad \bfphi \Big) } \dt = - 
	\intO{ \vr_0 \vc{w}_0 \cdot \bfphi (0, \cdot) },
\end{equation} 
for any test function $\bfphi$ in the class 
\begin{equation} \label{I3}
	\bfphi \in C^1_c([0,T) \times \Ov{\Omega}; R^d),\ \bfphi \cdot \vc{n}|_{\partial \Omega} = 0.
\end{equation}	
For the proofs from the previous sections to be adaptable to this boundary conditions, we have to impose a geometric restriction on the shape of the domain $\Omega$, specifically, $\Omega$ \emph{is not} rotationally symmetric with respect to some axis.
\end{Remark}

\begin{Remark}
We can also consider the case of general boundary conditions on a regular bounded domain $\Omega$, namely, 
\begin{equation} \label{G1}
	\vr \vu \cdot \vc{n}|_{\partial \Omega} = v_B,\ 
	\vr|_{\partial \Omega} = \vr_B,\ \Grad \vr \cdot \vc{n}|_{\partial \Omega} = D_N \vr_B,
\end{equation}
where $D_N$ is the Dirichlet-to-Neumann operator.\\
To accommodate \eqref{G1}, we consider a weaker formulation of \eqref{i2}, namely
\begin{equation} \label{G2}
	\int_0^T \intO{ \Big( \vr \vc{w} \cdot \partial_t \bfphi + 
		\vr \vc{w} \otimes \vu : \Grad \bfphi \Big) } \dt = - 
	\intO{ \vr_0 \vc{w}_0 \cdot \bfphi (0, \cdot) }
\end{equation} 
for any test function $\bfphi$ in the class 
\begin{equation} \label{G3}
	\bfphi \in C^1_c([0,T) \times {\Omega}; R^d).
\end{equation}

In the case of general boundary conditions, we can obtain the existence of infinitely many solutions for 
given data. The related energy inequality must be modified accordingly to discuss admissible solutions 
for certain data. 
\end{Remark}

\def\cprime{$'$} \def\ocirc#1{\ifmmode\setbox0=\hbox{$#1$}\dimen0=\ht0
	\advance\dimen0 by1pt\rlap{\hbox to\wd0{\hss\raise\dimen0
			\hbox{\hskip.2em$\scriptscriptstyle\circ$}\hss}}#1\else {\accent"17 #1}\fi}


\begin{thebibliography}{99}


\bibitem{AwKlar}
A.~Aw, A.~Klar, M.~Rascle, and T.~Materne.
\newblock Derivation of continuum traffic flow models from microscopic
  follow-the-leader models.
\newblock {\em SIAM J. Math. Anal.}, 63(1):259--278, 2002.

	
	\bibitem{AR}
	A. Aw and M. Rascle. 
	\newblock Resurrection of second order models of traffic flow.
	\newblock {\em SIAM J. Appl. Math.}, 60:916--938, 2000.
	
	
	\bibitem{BFHZ}
D. Breit, E. Feireisl, M. Hofmanova, and E. Zatorska.
\newblock Compressible Navier--Stokes system with transport noise.
\newblock {\em SIAM J. Math. Anal.}, 54(4), 937--972, 2022. 

\bibitem{BVY}
D. Bresch, A. Vasseur, C.Yu.
\newblock Global existence of entropy-weak solutions to the compressible Navier-Stokes equations with nonlinear density dependent viscosities. \newblock To appear in {\em J. Eur. Math. Soc.}, arXiv:1905.02701, 2019.
	
	\bibitem{BDDR}
	 F. Berthelin, P. Degond, M. Delitata, and M. Rascle. 
	 \newblock A Model for the Formation and Evolution of Traffic Jams. 
	 \newblock{\em Archive for Rational Mechanics and Analysis}, 187, 185--220, 2008. 
	 
\bibitem{BuVi}
T. Buckmaster and V. Vicol.
\newblock Nonuniqueness of weak solutions to the Navier-Stokes equation.
\newblock{\em Ann. of Math.} (2) 189(1): 101-144, 2019.



\bibitem{CVY}
R. M. Chen, A. F. Vasseur, C. Yu.
\newblock Global ill-posedness for a dense set of initial data to the Isentropic system of gas dynamics.
\newblock To appear in {\em Advances in Mathematics},	arXiv:2103.04905, 2021.


\bibitem{CGZ}
N. Chaudhuri, P. Gwiazda, E. Zatorska.
\newblock Analysis of the generalised Aw-Rascle model
\newblock arXiv:2202.04130, 2022.

\bibitem{CFGG}
F. A. Chiarello, J. Friedrich, P. Goatin, and S. G\"ottlich.
\newblock Micro-Macro Limit of a Nonlocal Generalized Aw-Rascle Type Model.
\newblock {\em SIAM J. Appl. Math.}, 80 (4), 1841--1861, 2020.

\bibitem{Ch}
E. Chiodaroli. 
\newblock A counterexample to well-posedness of entropy solutions to the compressible Euler system. 
\newblock {\em J. Hyperbolic Differ. Equ.}, 11(3):493–519, 2014.

% 	\bibitem{CaWrZa2}
% J.~A. Carrillo, A.~Wr{\'o}blewska-Kami{\'n}ska, and E.~Zatorska.
% \newblock Pressureless {E}uler with nonlocal interactions as a singular limit
%   of degenerate {N}avier-{S}tokes system.
% \newblock {\em J. Math. Anal. Appl.}, 492(1), 2020.


%  \bibitem{PAM}
% B.~Haspot.
% \newblock { From the highly compressible {N}avier--{S}tokes equations to fast
%   diffusion and porous media equations, existence of global weak solution for
%   the quasi-solutions.}
% \newblock {\em J. Math. Fluid Mech.}, 18(2):243--291, 2016.

\bibitem{DeSz}
C. De Lellis and L. Sz\'ekelyhidi, Jr. 
\newblock On admissibility criteria for weak solutions of the Euler equations. 
\newblock {\em Arch. Ration. Mech. Anal.}, 195(1):225–260, 2010.

\bibitem{Fei2016}
E.~Feireisl.
\newblock Weak solutions to problems involving inviscid fluids.
\newblock In {\em Mathematical Fluid Dynamics, Present and Future}, Volume 183
of {\em Springer Proceedings in Mathematics and Statistics}, pages 377--399.
Springer, New York, 2016.

\bibitem{FeKlKrMa}
E.~Feireisl, C.~Klingenberg, O.~Kreml, and S.~Markfelder.
\newblock On oscillatory solutions to the complete {E}uler system.
\newblock {\em J. Differential Equations}, {\bf 269}(2):1521--1543, 2020.

 
 \bibitem{HZ}
B.~Haspot and E.~Zatorska.
\newblock { From the highly compressible {N}avier-{S}tokes equations to the
  porous medium equation -- rate of convergence.}
\newblock {\em Discrete Contin. Dyn. Syst.}, 36(6):3107--3123, 2016.

\bibitem{HMV}
M. Herty, S. Moutari, G. Visconti. 
\newblock Macroscopic modeling of multilane motorways using a two-dimensional second-order model of traffic flow.
\newblock {\em SIAM J. Appl. Math.} Vol. 78, No. 4, pp. 2252--2278, 2018.


\bibitem{LM}
 A. Lefebvre-Lepot and B. Maury,  
\newblock Micro-Macro Modelling of an Array of Spheres Interacting Through Lubrication Forces.
\newblock {\em Advances in Mathematical Sciences and Applications},
 vol. 21, nb 2, 535--557, 2011.

\bibitem{LW}
M. J. Lighthill and J. B. Whitham. 
\newblock {On kinematic waves: I. Flow movement
in long rivers. II. A theory of traffic flow on long crowded roads.}
\newblock{\em Proc. Roy. Soc.}, A229:1749--1766, 1955.
  
% \bibitem{PLL} 
% P.-L. Lions. 
% \it Mathematical Topics in Fluid Mechanics, Vol 2: Compressible Models. 
% \rm Oxford University Press, \rm New York, 1998.

\bibitem{MPI}	
S. Mueller
\newblock {\em	Variational  models for microstructure and phase transitions} 
\newblock{ S. Hildebrandt, M. Struwe (Eds.) Calculus of Variation and Geometric Evolution Problem, Lecture Notes in Math., vol. 1713, Springer-Verlag, Berlin Heidelberg (1999)}


\bibitem{Payne}
H. J. Payne.
\newblock { FREFLO: A macroscopic simulation model of freeway traffic}.
\newblock {\em Transportation Research Record}, 722:68--75, 1979.

\bibitem{R}
P. I. Richards.
\newblock { Shock waves on the highway.}
\newblock {\em Operations Research}, 4:42--51, 1956.

% \bibitem{tad4}
% R.~Shvydkoy and E.~Tadmor.
% \newblock { Eulerian dynamics with a commutator forcing {III}. {F}ractional
%   diffusion of order $0<\alpha<1$.}
% \newblock {\em Phys. D}, 376-377:131 -- 137, 2018.

% \bibitem{VaYu}
% A.~F. Vasseur and C.~Yu.
% \newblock {Existence of global weak solutions for 3d degenerate compressible {N}avier--{S}tokes equations.}
% \newblock {\em Invent. Math.}, 206(3):935--974, 2016.


\bibitem{Whitham}
G. B. Whitham.
\newblock {\em Linear and nonlinear waves.}
\newblock Wiley, New York, 1974.



 \end{thebibliography}
\end{document}